\documentclass[psamsfonts]{amsart}

\usepackage{amssymb,amsfonts}
\usepackage[all,arc]{xy}
\usepackage{enumerate}
\usepackage{mathrsfs}
\usepackage{amsrefs}
\usepackage{microtype}

\newtheorem{thm}{Theorem}[section]

\theoremstyle{definition}
\newtheorem{defn}[thm]{Definition}

\newtheorem{exmp}[thm]{Example}

\theoremstyle{remark}

\DeclareMathOperator{\trop}{trop}
\newcommand{\val}{\mathfrak{v}}

\makeatletter
\let\c@equation\c@thm
\makeatother
\numberwithin{equation}{section}

\begin{document}

\title{Computing Tropical Varieties in Macaulay2}

\author{Carlos Am\'endola \and Kathl\'en Kohn \and Sara Lamboglia \and Diane Maclagan \and Ben Smith\and Jeff Sommars \and Paolo Tripoli \and Magdalena Zajaczkowska}

\email{carlos.amendola@tum.de}
\email{kathlenk@math.uio.no},
\email{lambogli@math.uni-frankfurt.de}
\email{D.Maclagan@warwick.ac.uk}
\email{benjamin.smith@qmul.ac.uk}
\email{sommars1@uic.edu}
\email{paolo.tripoli@nottingham.ac.uk}
\email{Magdalena.A.Zajaczkowska@gmail.com}
\date{}

\begin{abstract}
We introduce a package for doing tropical computations in
\verb|Macaulay2|. The package draws on the functionality of \verb|Gfan| and
\verb|Polymake| while making the process as simple as possible for the end
user. This provides a powerful and user friendly tool for computing
tropical varieties requiring little prerequisite knowledge.
\end{abstract}

\maketitle

\markleft{}

\section{Introduction}
The purpose of this package is to facilitate computations in tropical geometry using \verb|Macaulay2| \cite{M2}.  
At the moment the main 
computational tool for 
tropical geometry is the program
\verb|Gfan| \cite{Gfan} by Jensen. This computes the Gr\"obner fan of
an ideal $I$ and includes functions to compute only the subfan of the
Gr\"obner fan given by the tropicalization of the variety $V(I)$.  The
polyhedral geometry program \verb|Polymake| \cite{GJ00} also has some
tropical functionality that is not implemented in \verb|Gfan|.

The package \verb|gfanInterface| \cite{gfanInterface}, implemented in
\verb|Macaulay2|, allows the user to interface with \verb|Gfan| while
retaining the computational speed provided by \verb|Macaulay2| for
Gr\"obner basis computations.  The program \verb|Gfan| is distributed
with \verb|Macaulay2|.  A drawback of the \verb|gfanInterface| package
is that it requires good knowledge of the functions and conventions of
\verb|Gfan|.  The goal of the \verb|Tropical| package is to provide a
user friendly tool to do these computations in \verb|Macaulay2|
without requiring any knowledge of these conventions.  The package
includes different strategies for the same function depending on the
input, and calls functions from \verb|Gfan|, via \verb|gfanInterface|,
and \verb|Polymake|, as appropriate.  Moreover the package implements
some extra functionality not yet available in \verb|Gfan|, such as
computing multiplicities for tropical varieties of non-prime ideals
and allowing the user to swap between the min and max conventions.

\section{Mathematical Background}

We follow the conventions of  chapters two and three of \cite{M-S}.
Let $K$ be a field with valuation $\val$.
\begin{defn}
Let $f=\sum_{{\bf u} \in \mathbb{Z}^n} a_{\textbf u}x^{\textbf u}$ be a polynomial
in $S=K[x_1^{\pm 1},\ldots,x_n^{\pm 1}]$. The \textit{tropicalization} of $f$ is
the function $\trop(f):\mathbb R^n\to \mathbb R$ given by 
\[
\trop(f)(\textbf w)=\min\{\val (a_{\textbf u})+\textbf w\cdot \textbf u : a_{\textbf u}\neq 0\}.
\]
The \textit{tropical hypersurface} defined by $f$ is 
\[
\trop(V(f))=\{\textbf w\in \mathbb R^{n} :\ \text{the minimum in }\trop(f)(\textbf w)\ \text{is achieved at least twice}\}.
\]
Let $I$ be an ideal in $S$. The \textit{tropicalization}
of the variety $V(I)$ is
\[
\trop(V(I))=\bigcap_{f\in I}\trop(V(f)).
\]
\end{defn}

The same definitions can be formulated using $\max$ instead of $\min$.
The \verb|Tropical| package allows the user to choose their convention
when loading the package.

If the ideal $I$ is generated by $f_1,...,f_s$ then it is
not true in general that $\trop V(I)$ is the intersection
of the tropical hypersurfaces associated to the $f_i$s.
The intersection $\trop V(f_1)\cap ...\cap \trop V(f_s)$
is a \textit{tropical prevariety}.

\begin{defn}
Let $I=\langle f_1,...,f_s \rangle$ be an ideal in $S=K[x_1^{\pm 1},\ldots,x_n^{\pm 1}]$.
Then $f_1,...,f_s$ are a \textit{tropical basis} of $I$
if $\trop V(I)=\bigcap _{i=1}^s \trop V(f_i)$. 
\end{defn}

The tropical variety $\trop V(I)$ is a 
polyhedral complex
(\cite{M-S}*{Proposition 3.2.8}) contained in the Gr\"obner complex 
of the ideal $I$. If the valuation is trivial, the
tropical variety is a rational polyhedral fan and is a subfan of the Gr\"obner
fan of the ideal. Moreover we can associate to each maximal cell an integer
number called a \textit{multiplicity} such that a balancing condition
holds (see \cite{M-S}*{Definition 3.3.1}).  

The \verb|Tropical| package takes as input an ideal $I$ in a usual
(non-Laurent) polynomial ring $K[x_1,\dots,x_n]$.  The tropical
variety computed is the variety $V(J) \subset (K^*)^n$ of the ideal
$J=IK[x_1^{\pm 1},\dots,x_n^{\pm 1}]$.

\section{Examples}

In this section we give explicit examples in order to give
a short overview of the package. The computations are all
over the field $\mathbb Q$ of rational numbers with trivial
valuation, hence all tropical varieties are polyhedral fans.

\begin{exmp}
Consider the algebraic variety $X=V(I)\subset (\mathbb C^*)^2$, where
$I=\langle x+y+1 \rangle$. The tropicalization of this variety can be
computed using the function \verb|tropicalVariety(I)|. The package
outputs this as a \emph{tropical cycle}: a fan with a list of multiplicities
corresponding to integer weights on the maximal cones. We extract
information about the tropical cycle using associated functions. For
example \verb|rays| gives the generators of the rays as the columns of
a matrix.
\begin{center}
    \begin{tabular}{ |p{6cm}   | p{6cm} |}
    \hline
    \begin{verbatim}
i1 : needsPackage("Tropical",
       Configuration=>{
       "tropicalMax"=>false});
i2 : R=QQ[x,y];
i3 : I=ideal(x+y+1);
i4 : T=tropicalVariety I; 
o4 = T
o4 : TropicalCycle
i5 : rays T
o5 = | -1 1 0 |
     | -1 0 1 |
              2        3
o5 : Matrix ZZ  <--- ZZ
 \end{verbatim} & \begin{verbatim}
i6 : linealitySpace T
o6 = 0
              3
o6 : Matrix ZZ  <--- 0
i7 : maxCones T
o7 = {{0}, {1}, {2}}
o7 : List
i8 : multiplicities T
o8 = {1, 1, 1}
o8 : List
\end{verbatim}  \\ \hline

    \end{tabular}
\end{center}

The 
tropical variety $\trop V(I)$ is the standard tropical line in the plane: a $1-$dimensional fan
in $\mathbb R^2$ whose rays are $(-1,-1), (1,0)$, and $(0,1)$.
\end{exmp}

The function \verb|tropicalVariety| uses one of two different
 algorithms depending on the input ideal. If the ideal is prime, the
 tropical variety is connected through codimension one
 (\cite{BJSST}*{Theorem 3.1}) and the \verb|Gfan| commands
 \verb|gfan_tropicalstartingcone| and \verb|gfan_tropicaltraverse|,
 which implement the algorithm described in \cite{BJSST}, are used.
  However if the ideal is not prime, this algorithm might give the
 wrong answer.  The package then calls the more computationally
 expensive command \verb|gfan_tropicalbruteforce|, which computes the
 entire Gr\"obner fan.  The multiplicities are then computed
 separately.
The package does not require that
 the user knows these intricacies, but simply requires that they flag
 when the ideal is not prime.
\begin{center}
\begin{tabular}{ |p{6cm}   | p{6cm} |}
\hline
\begin{verbatim}
i9 : elapsedTime(
       tropicalVariety I);
     -- 0.088835 seconds elapsed
\end{verbatim} & \begin{verbatim}
i10 : elapsedTime(
        tropicalVariety(
          I,Prime=>false));
     -- 0.103651 seconds elapsed
\end{verbatim} \\ \hline
\end{tabular}
\end{center}

For most functions  \verb|Gfan| requires the input to be homogeneous.
 The \verb|Tropical| package will accept non-homogeneous
input, and do the pre- and post-processing to put it into a format
acceptable for \verb|Gfan|.   Small additions such as this help
decrease the prerequisite knowledge for the package.

\begin{exmp}
A tropical cycle is a fan with multiplicities attached to its maximal
cones; it need not be the tropicalization of an algebraic
variety. Therefore the package allows the user to create a tropical
cycle manually by defining a fan via its maximal cones and attaching
multiplicities to each of those cones. The following example shows how
we can construct $\trop V(I)$ manually.

\begin{center}
    \begin{tabular}{ |p{6.5cm}   | p{5.5cm} |}
    \hline
\begin{verbatim}
i11 : C1=posHull(matrix{{1},{0}});
i12 : C2=posHull(matrix{{0},{1}});
i13 : C3=posHull(matrix{{-1},{-1}});
i14 : F=fan({C1,C2,C3})
o14 = F
o14 : Fan
 \end{verbatim} & \begin{verbatim}
i15 : mult={1,1,1}
o15 = {1, 1, 1}
o15 : List
i16 : S=tropicalCycle(F,mult)
o16 = S
o16 : TropicalCycle
i17 : isBalanced S
o17 : true
\end{verbatim}  \\ \hline
    \end{tabular}
\end{center}
The \verb|tropicalCycle| command does not check that the resulting weighted fan is balanced.  To verify this we use the \verb|isBalanced| command.

\end{exmp}

\begin{exmp}
Consider the tropical hypersurfaces $\trop V(f)$ and $\trop V(g)$
cut out by the polynomials $f = x+y+z$ and $g=x^2+y^2+z^2$. Their
intersection cuts out a tropical prevariety. We would like to compute
whether this prevariety is equal to the tropical variety $\trop V(I)$
where $I = \langle f,g \rangle$.
\begin{center}
    \begin{tabular}{ |p{6cm}   | p{6cm} |}
    \hline
\begin{verbatim}
i18 : R=QQ[x,y,z];
i19 : f=x+y+z;
i20 : g=x^2+y^2+z^2;
i21 : l={f,g};
i22 : Tp=tropicalPrevariety l;
i23 : Tv=tropicalVariety ideal l;
 \end{verbatim} & \begin{verbatim}
i24 : isTropicalBasis l
o24 = false
i25 : dim Tp
o25 = 2
i26 : dim Tv
o26 = 1
\end{verbatim}  \\ \hline
    \end{tabular}
\end{center}
The polynomials $f,g$ are not a tropical basis for $I$ and therefore
the prevariety given by them is not equal to $\trop V(I)$.  We can see 
from our computation that the prevariety has a two-dimensional cone, while $\trop
V(I)$ is one-dimensional.
\end{exmp}

\begin{exmp}
For two curves $V(f)$ and $V(g)$ in $\mathbb{P}^2$, B\'{e}zout's
Theorem states that $|V(f) \cap V(g)|$ equals $ \text{deg}(f) \cdot
\text{deg}(g)$ counting multiplicities. The tropical analogue of
B\'{e}zout's Theorem 
states that the \emph{stable
intersection} of $\trop V(f)$ and $\trop V(g)$ is $\text{deg}(f) \cdot
  \text{deg}(g)$ points counting multiplicities. The following example
  shows how the package and the \verb|stableIntersection| function can
  be used to verify examples of tropical B\'{e}zout's Theorem.
\begin{center}
    \begin{tabular}{ |p{6.5cm}   | p{5.5cm} |}
    \hline
\begin{verbatim}
i27 : f=random(2,R);
i28 : g=random(1,R);
i29 : Tf=tropicalVariety ideal f;
i30 : Tg=tropicalVariety ideal g;
i31 : Tint=stableIntersection(Tf,Tg)
o31 = Tint
o31 : TropicalCycle
i32 : rays Tint
 \end{verbatim} & \begin{verbatim}
o32 = 0
               3        
o32 : Matrix ZZ  <--- 0
i33 : maxCones Tint
o33 = {{}}
o33 : List
i34 : multiplicities Tint
o34 = {2}
o34 : List
\end{verbatim}  \\ \hline
    \end{tabular}
\end{center}
The above code considers the stable intersection of a tropical line
and a plane quadric. The resulting tropical cycle is a single point, the
origin, with multiplicity two, verifying the claim of tropical
B\'{e}zout's theorem.

The function \verb|stableIntersection| has two different strategies
for computation depending on the software available to the user. If
the user has \verb|Polymake| version 3.2 or later installed, the
default strategy is to use \verb|atint| \cite{Atint}, a
\verb|Polymake| extension for tropical intersection theory by Simon
Hampe. If this is not available, the package instead uses \verb|Gfan|
to compute the stable intersection.
\end{exmp}

\section{Future Plans}

We plan for the \verb|Tropical| package to become the umbrella package
for all tropical computations in \verb|Macaulay2|.  This will include
implementing alternate strategies for some of the core commands as
algorithms improve, before they are included into \verb|Gfan| and
\verb|Polymake|.  

In addition there are still functions available in \verb|Gfan| and
\verb|Polymake| that are not yet available in the package.  We
particularly highlight the treatment of nontrivial valuations, which
is available in \verb|Gfan|, and the visualization of low-dimensional
tropical varieties, which is available in \verb|Polymake|.

\subsection*{Acknowledgements}

We thank Mike Stillman and Dan Grayson for their support
and help with the inner workings of \verb|Macaulay2|, Lars Kastner and
Georg Loho for their assistance with \verb|Polymake| and the
\verb|Polyhedra| package, Anders Jensen for helping with the
interactions with \verb|Gfan|, and Josephine Yu for helpful discussions
on \verb|Polymake| and \verb|Gfan|.

\end{document}